\documentclass[12pt]{article}
\usepackage{amsmath}
\usepackage{amsfonts}
\usepackage{amssymb}
\usepackage{amsthm}
\usepackage{hyperref}

\usepackage{graphicx}

\def\bee{\begin{equation}}
\def\eee{\end{equation}}

\begin{document}

\thispagestyle{empty}
\centerline{}
\bigskip
\bigskip
\bigskip
\bigskip
\bigskip
\bigskip
\centerline{\Large\bf Two  arguments that the nontrivial zeros}
\bigskip
\centerline{\Large\bf of the Riemann zeta function are irrational. II}
\bigskip
\bigskip

\begin{center}
{\large \sl Marek Wolf}\\*[5mm]
\begin{center}
Cardinal  Stefan  Wyszynski  University, Faculty  of Mathematics and Natural Sciences. College of Sciences,\\
ul. W{\'o}ycickiego 1/3,   PL-01-938 Warsaw,   Poland, e-mail:  m.wolf@uksw.edu.pl
\end{center}

\bigskip
\end{center}

\bigskip
\bigskip

\begin{center}
{\bf Abstract}\\
\bigskip
\begin{minipage}{12.8cm}
We extend the results of our previous  computer experiment performed on  the first 2600 nontrivial
zeros $\gamma_l$ of the Riemann zeta function calculated  with 1000 digits accuracy  to the set of 40000
first zeros given  with 40000 decimal digits accuracy.   We  calculated the geometrical means of the denominators of
continued fractions  expansions  of these zeros and for all cases we get values very  close to the
Khinchin's constant, which  suggests that $\gamma_l$ are irrational. Next we
have calculated the $n$-th square roots of the denominators $Q_n$ of the
convergents of the continued fractions obtaining values very  close to the Khinchin---L{\'e}vy
constant, again supporting the common opinion  that $\gamma_l$ are irrational.
\end{minipage}
\end{center}

\bigskip\bigskip

\section{Introduction}

The famous Riemann's  $\zeta(s)$ function \cite{Edwards,  Titchmarsh}  has trivial zeros at even negative integers:  $-2, -4, -6, \ldots$
and infinity of  nontrivial complex zeros $\rho_l=\beta_l+i\gamma_l, ~\rho_{-l}=\beta_l-i\gamma_l=\overline{\beta_l+i\gamma_l}=\overline{\rho_l}$
($l=1, 2, 3, \ldots$)  (here overline  denotes complex  conjugation:  $\overline{z}=\overline{x+iy}=x-iy$)   in the critical strip: $\beta_l\in (0,1),~\gamma_l \in \mathbb{R}$.
The Riemann Hypothesis (RH) asserts that $\beta_l=\frac{1}{2}$  for all $l$  --- i.e.
all zeros lie on the critical line $\Re(s)=\frac{1}{2}$.  Presently it is added  that
these nontrivial zeros are simple: $\zeta'(\rho_l)\neq 0$ --- many explicit
formulas of number theory contain $\zeta'(\rho_l)$ in the denominators.

There is little hope to obtain an analytical  formulas for the imaginary parts $\gamma_l$ of
the nontrivial zeros of $\zeta(s)$.  In the paper  \cite{Leclair2013}  LeClair obtained  almost  exact  equation
for $\gamma_l$  whose  solution  gives  approximate values of imaginary parts of  nontrivial   zeros in terms
of the Lambert function $W(x)$:

\[
x=W(x)e^{W(x)}
\]
see  e.g.  \cite[p. 111,  \S 4.13]{NIST-2010}
In the previous paper \cite{Wolf_2018}  we  exploited  two facts about the continued fractions: the
existence of the   Khinchin constant and  Khinchin--L{\'e}vy constant,
see e.g. \cite[\S 1.8]{Finch},  to support the irrationality of $\gamma_l$.  Let
\bee
r = [a_0(r); a_1(r), a_2(r), a_3(r), \ldots]=
a_0(r)+\cfrac{1}{a_1(r) + \cfrac{1}{a_2(r) + \cfrac{1}{a_3(r)+\ddots}}}
\eee
be the continued fraction expansion of the real number $r$, where $a_0(r)$ is an integer and all $a_k(r)$ with $k\geq 1$
are positive integers.  Khinchin has proved \cite{Khinchin}, see also \cite{Ryll-Nardzewski1951}, that  limits of geometrical
means of $a_n(r)$ are the same  for almost all real $r$:
\bee
\lim_{n\rightarrow \infty} \big(a_1(r) \ldots a_n(r)\big)^{\frac{1}{n}}=
\prod_{m=1}^\infty {\left\{ 1+\frac{1}{ m(m+2)}\right\}}^{\log_2 m} \equiv K_0 \approx 2.685452001\dots   ~~.
\label{Khinchin}
\eee
The Lebesgue measure of all exceptions is zero  and include {\it rational numbers}, quadratic irrationals and
some irrational numbers too, like for example the
Euler constant $e=2.7182818285\ldots$ for which the limit (\ref{Khinchin}) is
infinity.

The constant $K_0$ is called the Khinchin constant, see e.g.  \cite[\S 1.8]{Finch}.
If the quantities
\bee
K(r; n)=\big(a_1(r) a_2(r) \ldots a_n(r)\big)^{\frac{1}{n}}
\label{Kny}
\eee
for a given number $r$ are close to $K_0$ we can regard it as an indication  that $r$ is irrational.
This  is  the idea behind  our papers:  previous \cite{Wolf_2018}  and  the present one.

Let the rational $P_n/Q_n$ be the $n$-th  partial convergent of the continued fraction:
\bee
\frac{P_n}{Q_n}=[a_0; a_1, a_2, a_3, \ldots, a_n].
\eee
For almost all real numbers $r$ the denominators of the  finite continued fraction
approximations fulfill:
\bee
\lim_{n \rightarrow \infty} \big(Q_n(r)\big)^{1/n} = e^{\pi^2/12\ln2} \equiv L_0 = 3.275822918721811\ldots
\eee
where $L_0$ is called the  Khinchin---L{\'e}vy's constant  \cite[\S 1.8]{Finch}.
Again the set of exceptions to  the above limit is  of the Lebesgue measure zero and it includes rational numbers,
quadratic irrational etc.

\section{Computer experiments with  zeta zeros}

Some  time ago we learned that G. Beliakov  and Y. Matiyasevich    \cite{Beliakov_Matiyasevich_2013}  calculated first  40000
nontrivial zeros of $\zeta(s)$  with  40000  digits accuracy and made them  publicly available at \cite{zera_40000}.
The method  used    during this high-precision numerical calculations is described in  \cite{Beliakov_Matiyasevich_2015}.  
Thus we  were able to repeat  computer  experiments   from \cite{Wolf_2018}  on much larger set of zeros of $\zeta(s)$  given  with
much  higher  number of digits.

In the computer experiments we used the PARI \cite{PARI2}  which has built in function  \verb"contfrac"$(r,\{nmax\})$.
This function  creates the row vector ${\bf a}(r)$ whose
components are the denominators $a_n(r)$  of the continued fraction
expansion of $r$, i.e.  ${\bf a}=[a_0(r); a_1(r), \dots,a_n(r)]$ means that
\bee
r \approx
a_0(r)+\cfrac{1}{a_1(r) + \cfrac{1}{a_2(r) + \cfrac{1}{\ddots\cfrac{1}{a_n(r)}}}}.
\eee
The parameter $nmax$ limits the number of terms $a_{nmax}(r)$; if it is omitted
the expansion stops with a declared precision  of representation of real numbers at the last
significant partial quotient.    With the precision set  to 90000 digits we  expanded  each $\gamma_l$, $l=1, 2, \ldots 40000$
with  40000  accurate  decimal   digits   value  into  its   continued fractions
\bee
\gamma_l \doteq [a_0(l); a_1(l), a_2(l), a_3(l), \ldots, a_{n(l)}(l)]\equiv {\bf a}(l)
\label{g-cfr}
\eee
(here $\doteq$   denotes approximate equality)   without specifying the  parameter $nmax$, thus the length $n(l)$ of
the vector ${\bf a}(l)$  depended on $\gamma_l$ and it turns out that the number of denominators was contained between 77000 and  78000.
The value of the product  $a_1 a_2 \ldots a_{n(l)} $  was typically of the
order $10^{33000}---10^{33500}$.   Next for each $l$ we have calculated the geometrical means:
\bee
K_l(n(l))= \left( \prod_{k=1}^{n(l)} a_k(l) \right)^{1/n(l)}.
\eee
The results are presented in the Fig.1. Values of $K_l(n(l))$ are scattered around the
red line representing $K_0$  and are contained in the interval $(K_0-0.06, K_0+0.06)$.  For the set of  zeros
reported in the  paper  \cite{Wolf_2018}  values of  $K_l(n(l))$   were   contained in the interval $(K_0-0.3, K_0+0.3)$.
The Fig.1  is in some sense misleading,  because  there are cases  when the  difference $K_l(m)-K_0$ changes  sign
for earlier  $m$.
We are not able to repeat calculations  presented in  the Fig. 2 in \cite{Wolf_2018}  showing the number of sign changes
of $K_l(m)-K_0$:  it would take a few years of CPU time as  there are over 15 times more zeros  with 40 times more digits,
thus assuming linear  complexity  of the problem  it will take over 600  more time than in the previous experiment, which took
2 CPU days.   In the Fig. 2 we present plots of  $K_l(m)$ as a function of $m$  for a few zeros $\gamma_l$  with  sign changes of
$K_l(m)-K_0$.   We encountered also zeros $\gamma_l$   without  sign change of  $K_l(m)-K_0$, some of them are plotted in the Fig.3.
The plots of the difference  $K_0-K_l(m)$  without  sign  changes seem to  follow the power--like  dependence
$|K_l(m)-K_0|\sim m^{-\alpha_l}$,  where the parameters $\alpha_l$   very weakly depend on the zero number $l$  and
are close to $0.9$.

\begin{figure}[h]
\begin{center}
\includegraphics[width=0.8\textwidth, angle=0]{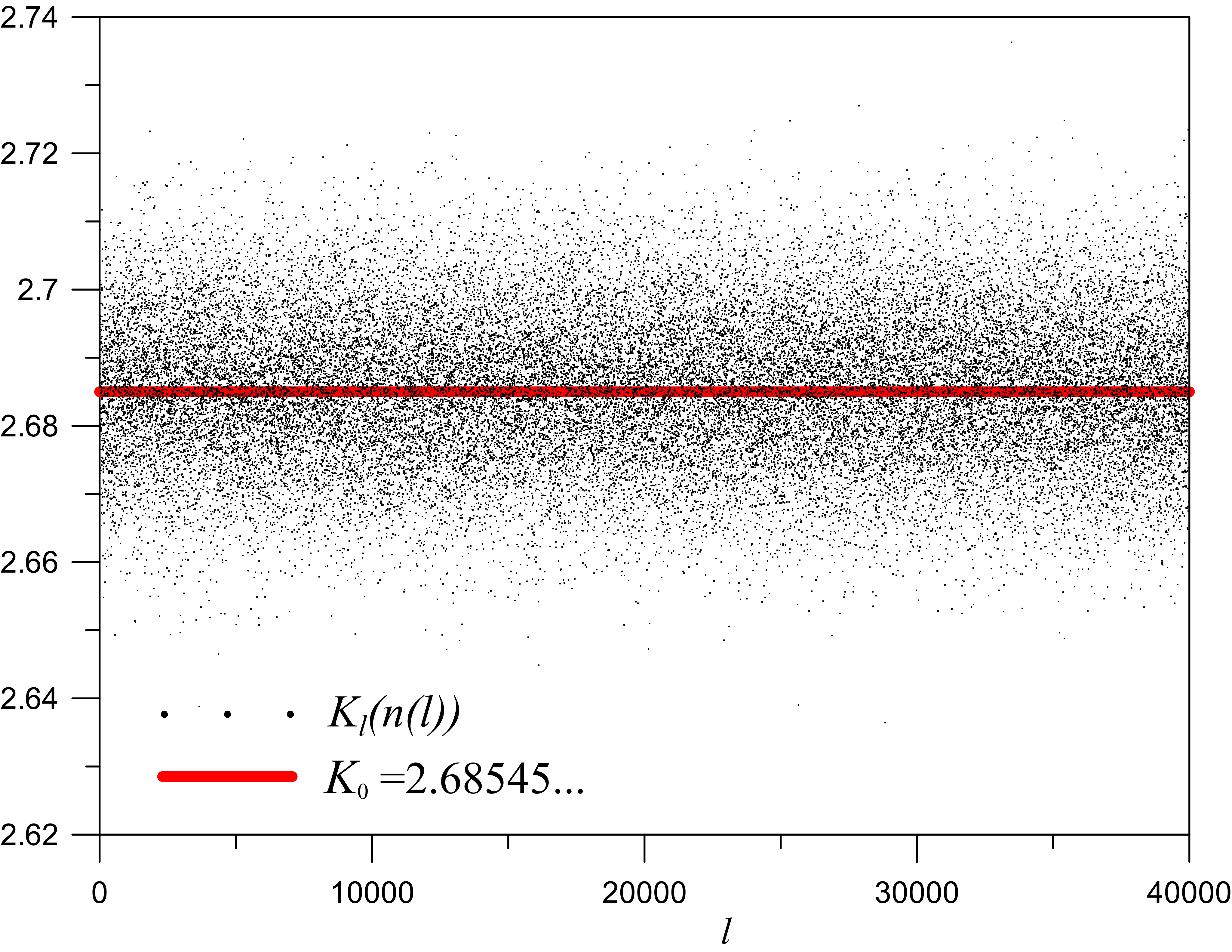} \\
\vspace{-1.5cm}
\vspace{1.7cm}Fig.1  The plot of $K_l(n(l))$ for $l=1,2,3, \ldots, 40000$.  There are 2976 points  closer  to $K_0$ than $0.001$
and 292 points  closer  to $K_0$ than $0.0001$.
The largest  value of  $|K_0-K_l(n(l))|$ is $5.08\times 10^{-2}$ and it occurred for the zero number $l=33473$,
the smallest  value of  $|K_0-K_l(n(l))|$ is  $9.2\times 10^{-8}$    and it occurred for the zero number $l=17408$.   \\
\end{center}
\end{figure}


\begin{figure}[h]
\begin{center}
\includegraphics[width=0.8\textwidth, angle=0]{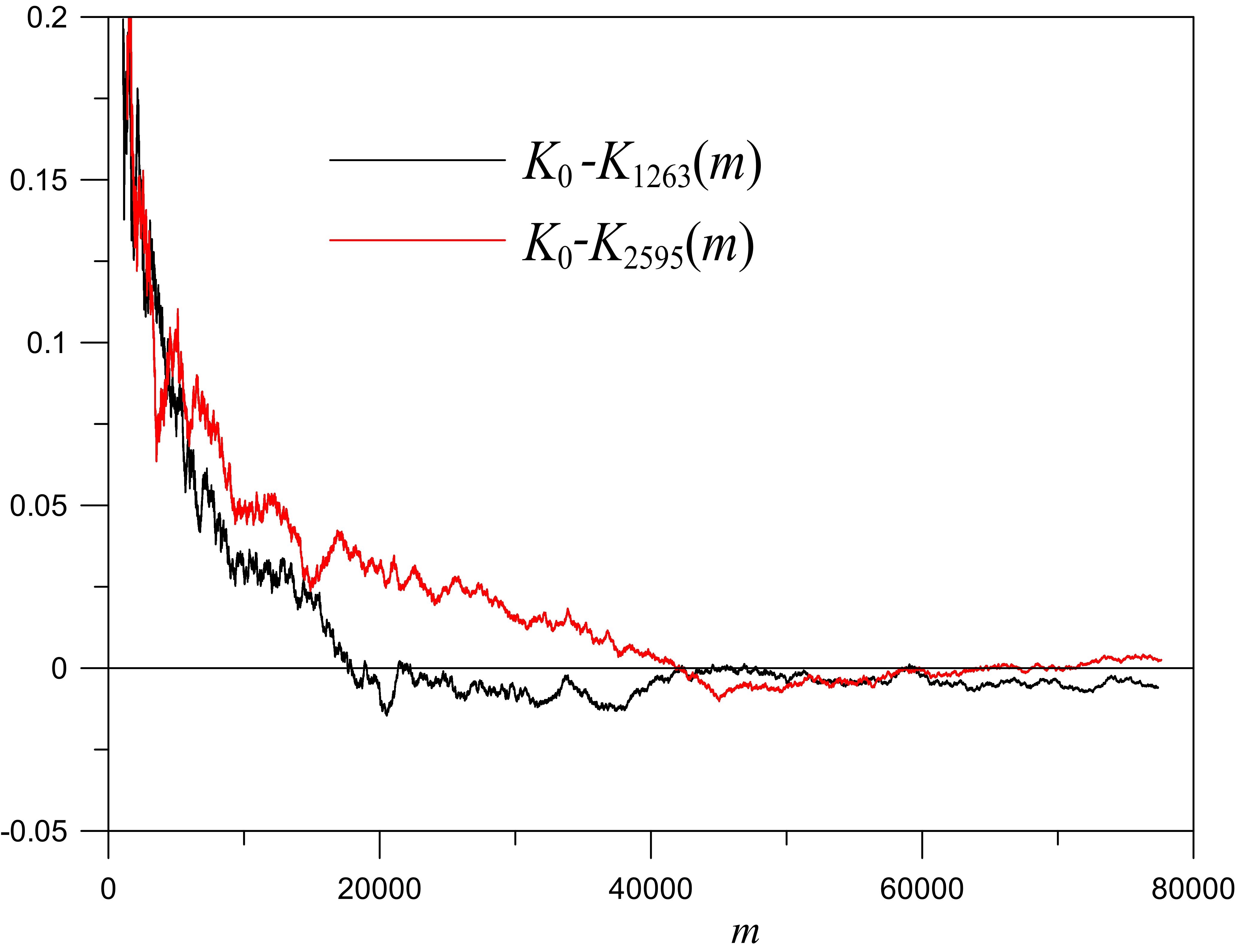} \\
\vspace{0.0cm}
Fig.2  The plot of the  difference $K_0-K_l(m)$ for $l=1263$ and $l=2595$.  There are 267  sign changes
for $\gamma_{1263}$   and 218 sign changes for $\gamma_{2595}$ \\
\end{center}
\end{figure}

We obtained from  G. Beliakov the zero $\gamma_{14299}$ with 70000  digits accuracy.  In the Fig.4 we present
the plot  of  $K_{14299}(m)-K_0$  vs $m$.   It took  52 h on 3.9 GHz CPU  to get data for this plot.  There are 322  sign changes
of  $K_0-K_{14299}(m)$ present on Fig.4.

Let the rational $P_{n(l)}(\gamma_l)/Q_{n(l)}(\gamma_l)$ be the $n$-th
partial convergent of the continued fractions (\ref{g-cfr}):
\bee
\frac{P_{n(l)}(\gamma_l)}{Q_{n(l)}(\gamma_l)}={\bf a}(l) \doteq\gamma_l .
\eee
For each zero $\gamma_l$ using PARI function \verb+contfracpnqn(a)+  we  calculated the partial convergents
$P_{n(l)}(\gamma_l)/Q_{n(l)}(\gamma_l)$.
Next from these denominators $Q_{n(l)}(\gamma_l)$ we have calculated the quantities
$L_l(n(l))$:
\bee
L_l(n(l))= \left(Q_{n(l)}\right)^{1/n(l)}, ~~~l=1, 2, \ldots , 40000
\eee
The obtained values of $L_l(n(l))$  are presented in the Fig.5.
These  values scatter around the red line representing the Khinchin---L{\'e}vy's
constant $L_0$  and are contained in the interval $(L_0-0.05, L_0+0.05)$, while in the previous paper  \cite{Wolf_2018}
this interval was $(L_0-0.36, L_0+0.36)$.   Again this plot is  somehow misleading
because there are zeros $\gamma(l)$  for which there appear  sign changes of  $L_0-L_l(m)$.
In the Fig. 4 we present the plot of $L_0-L_{14299}(m)$  for the zero number 14299 which is  known with 70000 digits
accuracy.  This plot is   practically identical  with  the plot of $K_0-K_{14299}(m)$. There are 229  sign changes
of  $L_0-L_{14299}(m)$ present on Fig.4.  These  sign changes  of  $K_0-K_{14299}(m)$   and  of  $L_0-L_{14299}(m)$  appears
almost  at the same arguments  $m$, see Fig.6.

{\bf  Acknowledgement:}  I thank  Gleb Beliakov for sending me the zero 14299  with 70000 digits accuracy.

\begin{figure}[h]
\begin{center}
\includegraphics[width=0.8\textwidth, angle=0]{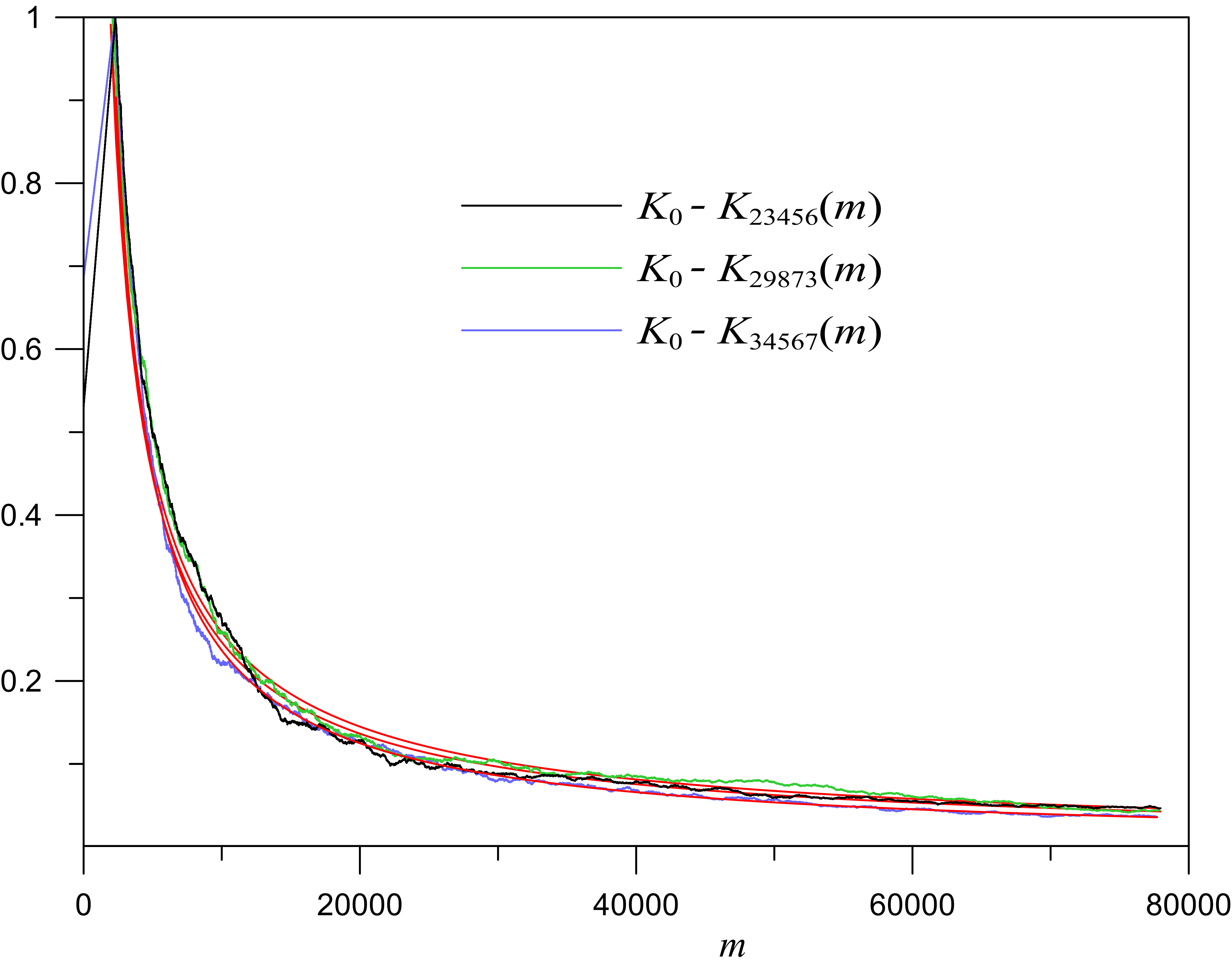} \\
\vspace{-0.0cm}
Fig.3  The plots of the  difference $K_0-K_l(m)$ for $l=23456, ~29873, ~34567$.  In red are the
fits to the power--like dependence plotted for these three  zeros.  These  fits  are  represented  by $m^{-\alpha}$
with parameter $\alpha$  almost the  same for all three zeros and contained in the interval $(0.85, 0.92)$.\\
\end{center}
\end{figure}

\begin{figure}[h]
\vspace{-0.3cm}
\begin{center}
\includegraphics[width=0.8\textwidth, angle=0]{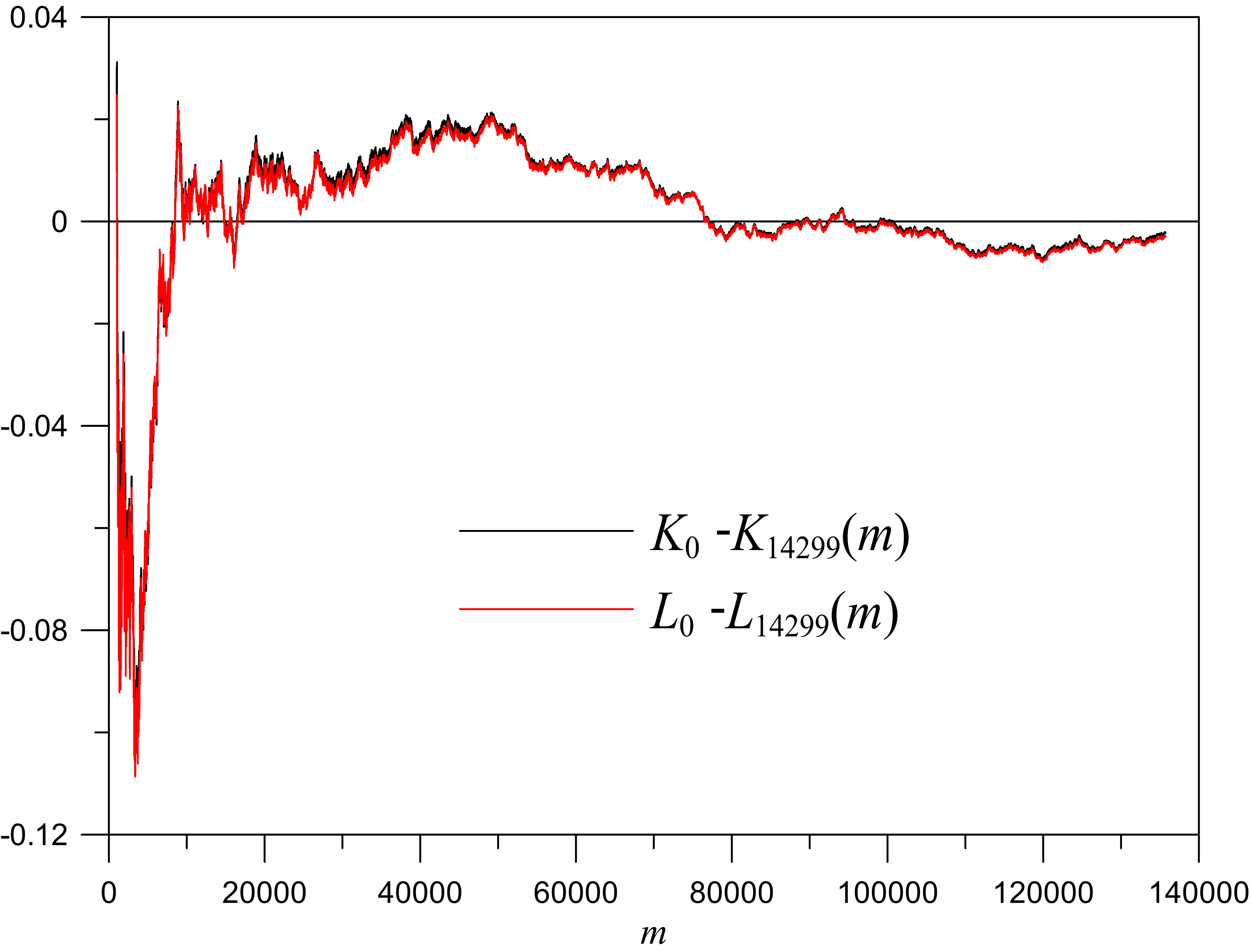} \\
\vspace{0.0cm}
Fig.4  The plots of the  differences  $K_0-K_{14299}(m)$  (black) and  $L_0-l_{14299}(m)$ (red).  The  number of denominators
$a(n)$  was  $135721$.\\
\end{center}
\end{figure}


\begin{figure}[h]
\vspace{-0.3cm}
\begin{center}
\includegraphics[width=0.8\textwidth, angle=0]{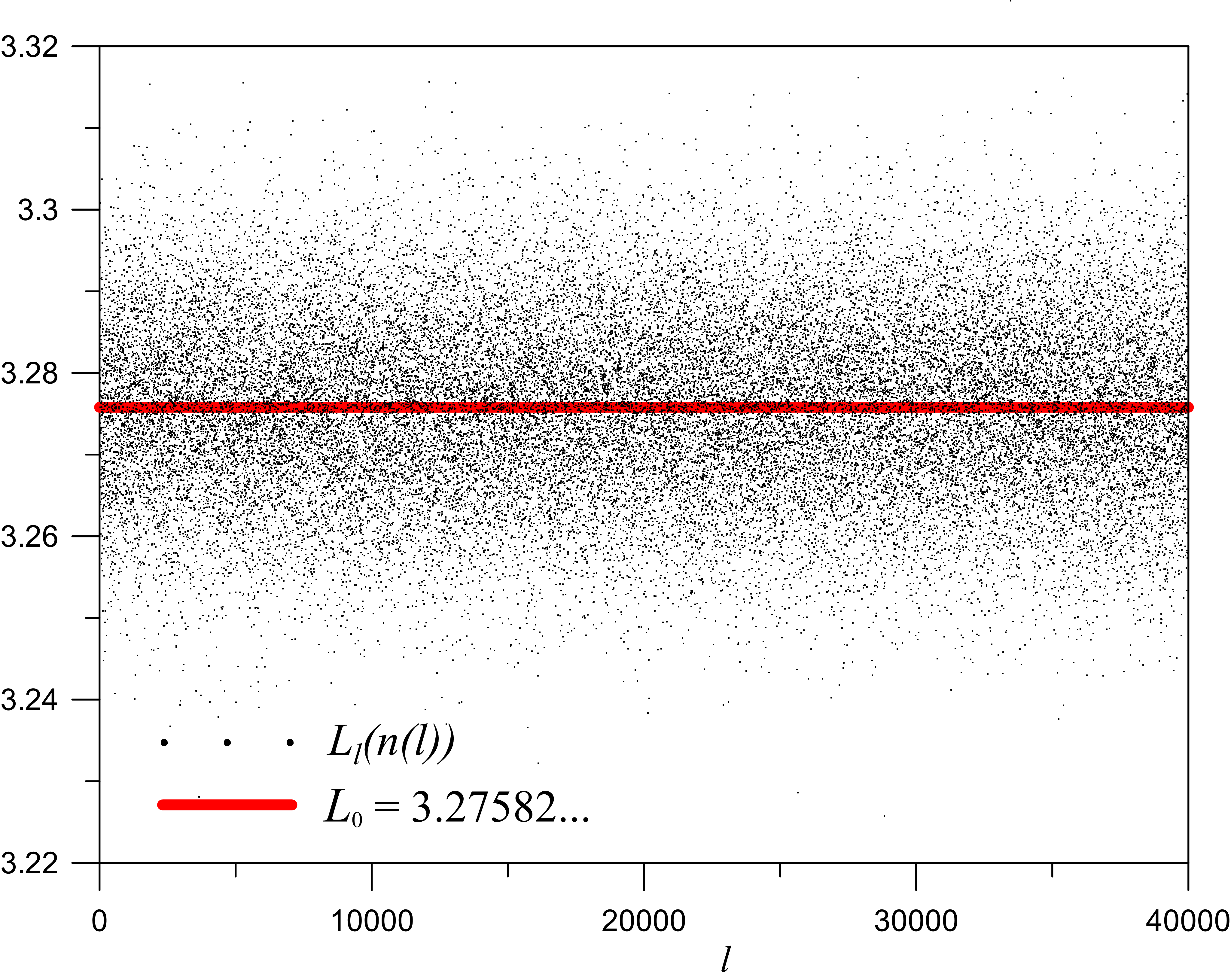} \\
\vspace{0.0cm}
Fig.5  The plot of $L_l(n(l))$ for $l=1,2,3, \ldots, 40000$.  There are 2906 points  closer  to $L_0$ than $0.001$ and 275  zeros
closer  to $L_0$ than $0.0001$  \\
\end{center}
\end{figure}

\begin{figure}[h]
\vspace{-0.3cm}
\begin{center}
\includegraphics[width=0.8\textwidth, angle=0]{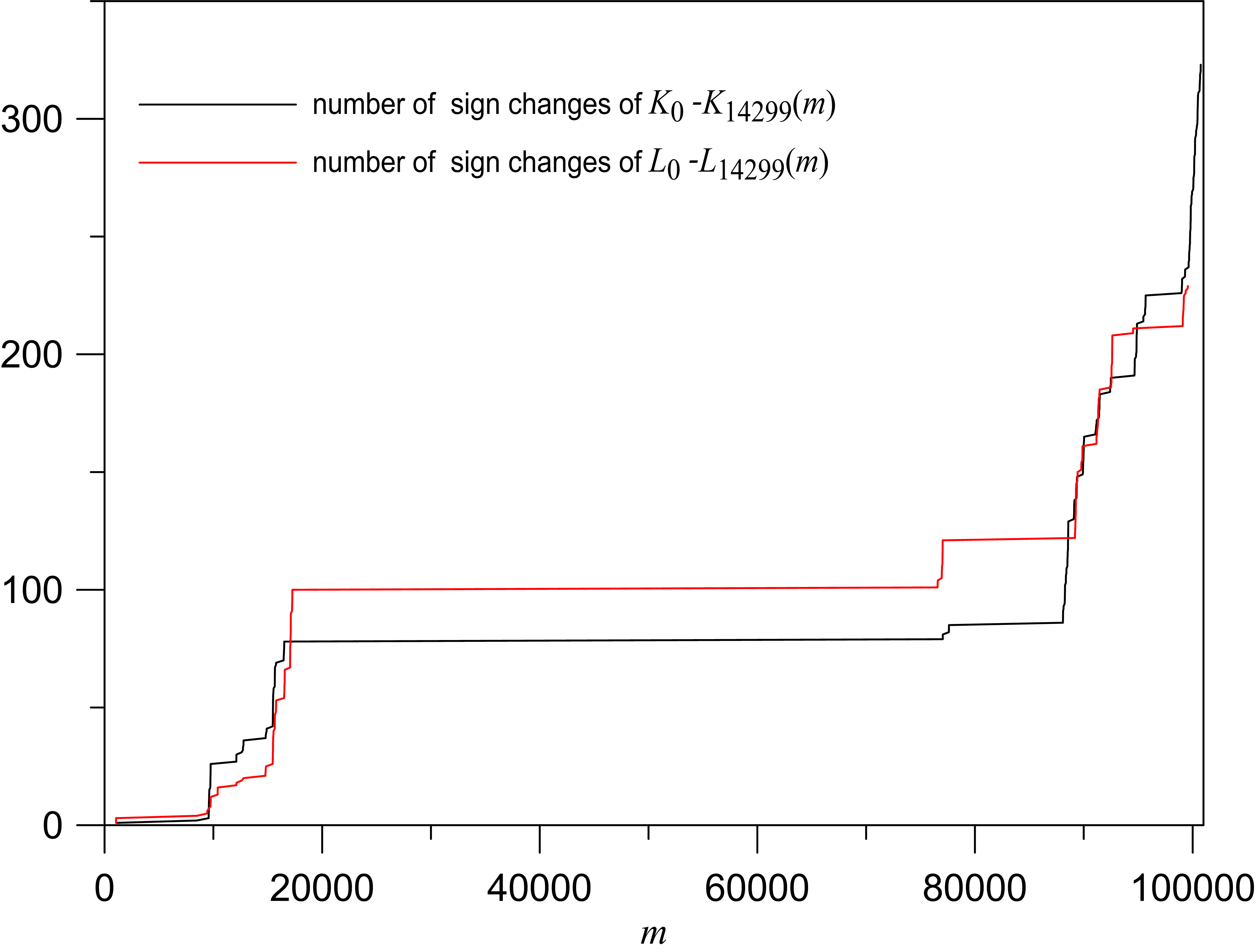} \\
\vspace{0.0cm}
Fig.6  The plot of the  number of sign changes of   difference $K_0-K_{14299}(m)$ and $L_0-L_{14299}(m)$ as  the function of $m$.\\
\end{center}
\end{figure}

\end{document}